\documentclass[conference]{IEEEtran}
\IEEEoverridecommandlockouts

\usepackage{graphics}
\usepackage{graphicx} 
\usepackage{rotating}
\usepackage{dcolumn}
\usepackage{hyperref}
\usepackage{longtable}
\usepackage{multirow}
\usepackage{amsmath}
\usepackage{amssymb}
\usepackage{xspace}
\usepackage{textcase}
\usepackage{widetext}
\usepackage{float}
\usepackage{amsthm}
\usepackage{multicol}
\usepackage{standalone}
\usepackage{tikz}
\usepackage{tikzscale}
\usepackage{pgfplots}
\pgfplotsset{compat=1.18}
\usetikzlibrary{trees}
\usetikzlibrary{arrows.meta}
\usetikzlibrary{positioning}
\usetikzlibrary{arrows, math, angles, quotes, shapes.geometric, calc}
\usepackage{commath}
\usepackage{xcolor}
\usepackage{hyperref}
\usepackage{subfiles}

\def\BibTeX{{\rm B\kern-.05em{\sc i\kern-.025em b}\kern-.08em
    T\kern-.1667em\lower.7ex\hbox{E}\kern-.125emX}}

\begin{document}

\title{Multi-Agent, Multi-Scale Systems with the Koopman Operator\
\thanks{This work was supported by the DOE Wind Energy Technologies Office.}
}

\author{\IEEEauthorblockN{1\textsuperscript{st} Craig Bakker}
\IEEEauthorblockA{\textit{National Security Directorate} \\
\textit{Pacific Northwest National Laboratory}\\
Richland, Washington \\
craig.bakker@pnnl.gov}
}

\maketitle

\begin{abstract}
The Koopman Operator (KO) takes nonlinear state dynamics and ``lifts'' those dynamics to an infinite-dimensional functional space of observables in which those dynamics are linear.  Computational applications typically use a finite-dimensional approximation to the KO.  The KO can also be applied to controlled dynamical systems, and the linearity of the KO then facilitates analysis and control calculations.  In principle, the potential benefits provided by the KO, and the way that it facilitates the use of game theory via its linearity, would suggest it as a powerful approach for dealing with multi-agent control problems.  In practice, though, there has not been much work in this space: most multi-agent KO work has treated those agents as different components of a single system rather than as distinct decision-making entities.  This paper develops a KO formulation for multi-agent systems that structures the interactions between decision-making agents and extends this formulation to systems in which the agents have hierarchical control structures and time scale separated dynamics.  We solve the multi-agent control problem in both cases as both a centralized optimization and as a general-sum game theory problem.  The comparison of the two sets of optimality conditions defining the control solutions illustrates how coupling between agents can create differences between the social optimum and the Nash equilibrium.
\end{abstract}


\section{Introduction}

The Koopman Operator (KO) takes nonlinear state dynamics and ``lifts'' those dynamics to an infinite-dimensional functional space of observables in which those dynamics are linear \cite{budisic12jsr}.  For a dynamical system $\dot{x} = f \left(x\right)$, the KO $\mathcal{K}$ advances an observable function $g \left(x\right)$ forward in time: $g \left(x_{t+1}\right) = \mathcal{K} g \left(x_t\right)$.  Computational applications typically use a finite-dimensional approximation to the KO: a matrix $K$, on a set of basis functions $\psi \left(x\right)$, producing $\psi \left(x_{t+1}\right) = K \psi \left(x_t\right)$ \cite{bakker2020koopman}.  The KO can also be applied to controlled dynamical systems, and the linearity of the KO then facilitates the calculation of optimal control policies using, e.g., Model Predictive Control (MPC) \cite{korda2018linear} or analyses such as observability quantification \cite{surana2016linear}.  When combined with linear constraints and a linear or quadratic cost function, the Koopman MPC is a linear or quadratic program, respectively, which makes it computationally efficient to solve.

In principle, this linearity would make the KO very useful for handling multi-agent dynamical systems using game theory (as a natural extension of single-agent control).  In practice, though, there has not been much work in this space.  Most multi-agent KO work has, for the purposes of control, treated those agents as different components of a single system (e.g., \cite{li2019learning,zhao2025koopman}).  Such ``agents'' are essentially different disciplines in a multi-disciplinary system, as described by Allison and Herber \cite{allison2014special} in the context of multidisiplinary design optimization: they are subsystems that have their own dynamics and interact with other subsystems.  This kind of multidiscplinary approach facilitates the incorporation of discipline structure as a kind of domain knowledge that can improve model performance, reduce data requirements, and/or enable useful analyses (e.g., symmetry as domain knowledge in \cite{sinha2020koopman}).  Multi-disciplinary KO models thus can to make a contribution to the broader field of domain-aware or physics-informed machine learning \cite{alber2019integrating}.

This is valuable in itself, but to treat an agents as \emph{agents}, not simply system components, implies \emph{agency} -- i.e., decision-making processes.  The consideration of interacting decision-making agents leads naturally to game theory, and there is a small amount of such game theoretic KO work in the literature.  Zhao and Zhu use the KO in a leader-follower game \cite{zhao2023stackelberg}, but they transform the game into a constrained optimization and only then use the KO to model the dynamics of the optimality conditions; using optimality conditions in a general nonlinear context raises some questions about the methodology.  Oster et al.~\cite{oster2023multi} combine the linearity of the KO with strong duality to convert a tri-level optimization into a bi-level optimization that they then solve with existing techniques from the multi-level optimization literature.  Bakker et al. \cite{bakker2025operator} use the KO in two fundamentally different ways to solve a zero-sum differential game: alternating optimizations using the resolvent of the Koopman operator and a Mixed Complementarity Problem (MCP), derived from the agents' optimality conditions, solved using the PATH algorithm \cite{dirkse1995path}.  Both the MCP and duality approaches rely on the linearity of the KO.  This would seem to motivate further use of KO-based control for multi-agent control problems, but to the best of our knowledge, only the examples listed above apply game theory and the KO together to dynamic games.

This paper presents structured KO formulations for multi-agent systems.  We consider agents, not just disciplines, but much of the paper can be applied, \textit{mutatis mutandis} to multi-disciplinary or multi-component systems (as defined above).  The structured formulations enable us to quantify the interactions between agents and how those interactions contribute to, e.g., overall system stability.  We then extend our parallel work on systems with hierarchical control and time scale separated dynamics \cite{bakker2025time} to multi-agent systems in which the agents have hierarchical control structures and time scale separated dynamics.  With this, we solve the multi-agent control problem both as a centralized optimization (i.e., one central controller determining the agents' control policies) and as a general-sum game theory problem (i.e., each agent calculating its own control policy).  Comparing the two sets of optimality conditions defining the control solutions illustrates how coupling between agents can create a difference between the social optimum and the Nash equilibrium (henceforth referred to as the optimum and equilibrium solutions, respectively).

\section{Koopman Model Description}

\subsection{Multi-Agent Systems}

Let us define a multi-agent system, with running costs, as

\begin{align}
		\dot{x}_i &= f_i \left(x_i,x_{-i},u_i\right) \\
		J_i &= \int \limits_0^T L_i \left(x_i,x_{-i}, u_i\right) dt
\end{align}

\noindent where $x_i$ is the state of agent $i$, $x_{-i}$ is the set of all $x_j, j \neq i$, $u_i$ is the control input to agent $i$, and $J_i$ is the objective of agent $i$.  Adding a terminal cost is straightforward but lengthens the derivations, so we omit it here.  We can then define a discrete-time Koopman formulation as

\begin{gather}
		\psi^x_i \left(x_{i,t+1}\right) = K^{xx}_{ii} \left(\psi^x_i \left(x_{i,t}\right) - \sum_{j \neq i} K^{xx}_{ij} \psi^x_j \left(x_{j,t}\right) \right. \nonumber \\
		\left. - K^{xu}_i \psi^u_i \left(u_{i,t}\right) \right) + \sum_{j \neq i} K^{xx}_{ij} \psi^x_j \left(x_{j,t}\right) + K^{xu}_i \psi^u_i \left(u_{i,t}\right) 
		\label{eq:discp dyn} \\
		%
		\boldsymbol \psi^x \left(\mathbf{x}_{t+1}\right) = K_{xx,comb} \boldsymbol \psi^x \left(\mathbf{x}_t\right)+ K_{xu,comb} \boldsymbol \psi^u \left(\mathbf{u}_t\right)
\end{gather}

\noindent where $\boldsymbol \psi \left(\mathbf{x}_t\right)$ is the stacked vector of all $\psi_i \left(x_{i,t}\right)$ vectors and $\boldsymbol \psi_u \left(\mathbf{u}_t\right)$ is the similarly stacked vector of all $\psi_{u,i} \left(u_{i,t}\right)$ vectors.  This multi-agent formulation could be applied to a multi-disciplinary or multi-component system simply by removing the control terms.  Here and in the rest of the paper, we use and modify the stability-assuming (or enforcing) KO formulation of King et al. \cite{king2021solving}.  The formulation in (\ref{eq:discp dyn}) separates out the agents such that each agent has its own set of observables $\psi_i$ and that the other agents' impacts on its dynamics are treated as inputs to the (agent-specific) system; the same goes for the control inputs associated with each agent.  We also define a quadratic discrete-time cost approximation $\hat{L}_k$ for agent $k$:

\begin{gather}
		\sum_t \left[\sum_{i,j} \left[\psi^i_{k,t}\right]^T Q_{k,ij} \psi^j_{k,t} + \sum_i c_{k,i}^T \psi^i_{k,t} + c_{k,0} \right]
		\label{eq:cost approx}
\end{gather}

\noindent where $i,j \in \left\{u,x\right\}$.  In principle, there are two ways that we could solve the multi-agent optimal control problem.  The first combines the agents into a single system (summing their objectives) and does a single optimization to solve for a social optimum.  The second way treats each agent as its own decision-making entity and solves for a Nash equilibrium.  This can be done by taking the optimality conditions of each agent's optimization (objective, dynamics, and constraints) and combining them together into a single MCP \cite{ruiz2014tutorial}.  If the agents' objectives are additive (as described here), the optimality conditions for the two solution methods are almost identical.  The one difference relates to the cross-agent coupling term in (\ref{eq:discp dyn}).  The equilibrium optimality condition in question is

\begin{align}
		0 &= \frac{\partial \hat{L}_k}{\partial \psi^x_{i,t}} - \lambda_{i,t-1} + \lambda_{i,t}^T K^{xx}_{ii} + \mu^x_{max,i,t} - \mu^x_{min,i,t}
\end{align}

\noindent where $\lambda_{i,t}$ is the dual variable associated with (\ref{eq:discp dyn}), $\mu^x_{max,i,t}$ is the dual variable associated with any upper bounds placed on $\psi^x_{i,t}$, and $\mu^x_{min,i,t}$ is the dual variable associated with any lower bounds placed on $\psi^x_{i,t}$.  The optimum optimality condition is

\begin{align}
		0 &= \frac{\partial \hat{L}_k}{\partial \psi^x_{i,t}} - \lambda_{i,t-1} + \lambda_{i,t}^T K^{xx}_{ii} + \sum_{j \neq i} \lambda_{j,t}^T \left(I - K^{xx}_{jj}\right) K^{xx}_{ji} \nonumber \\
		&+ \mu^x_{max,i,t} - \mu^x_{min,i,t}
\end{align}

\noindent so $\sum_{j \neq i} \lambda_{j,t}^T \left(I - K^{xx}_{jj}\right) K^{xx}_{ji}$ connects the solution across agents.  The equation does not indicate precisely how the optimum and equilibrium solutions will differ, but we can expect a difference.

\subsection{Multi-Agent Systems with Hierarchical Control and Time Scale Separation}

Let us now consider a multi-agent system where each agent has hierarchical control with time scale separated dynamics.  For agent $i$, we have a running cost and dynamics

\begin{gather}
		J_i = \int \limits_0^T L_i \left(w_i,y_i,u_i,x_i,x_{-i}\right) dt \\
		\begin{array}{c}\dot{x}_i = f_i \left(x_i,y_i,w_i,x_{-i}\right) \\
		\dot{y}_i = \frac{1}{\epsilon} g \left(x_i,y_i,w_i\right) \\
		\dot{w}_i = \frac{1}{\epsilon} h \left(x_i,y_i,w_i,u_i\right) \end{array} \xrightarrow{\epsilon \rightarrow 0} \begin{array}{c} \dot{x}_i = f \left(x_i,y_i,w_i,x_{-i}\right) \\
		0 = g \left(x_i,y_i,w_i\right) \\
		0 = h \left(x_i,y_i,w_i,u_i\right) \end{array} \\
		0 = g \left(x_i,y^*_i \left(x_i,u_i\right),w^*_i \left(x_i,u_i\right)\right) \\
		0 = h \left(x_i,y^*_i \left(x_i,u_i\right),w^*_i \left(x_i,u_i\right),u_i\right) 
\end{gather}

\noindent where $x_i$ is the slow scale state of agent $i$, $y_i$ is the fast scale state of agent $i$, $w_i$ is the (fast scale) actuators of agent $i$, $J_i$ is the running cost for agent $i$, and and $x_{-i}$ is the set of all $x_j, j \neq i$; see Bakker \cite{bakker2025time} for the prior single-agent version of such a system.  The Koopman model in (\ref{eq:x dyn})-(\ref{eq:w dyn}) takes the previously developed Koopman model for systems with hierarchical control and time scale separation \cite{bakker2025time} and combines it with the multi-agent model of the previous section; $\Delta t = m \tau$, $\tau \propto \epsilon$ captures the separation of scales.

\onecolumn
\par\noindent\rule{\dimexpr(0.5\textwidth-0.5\columnsep-0.4pt)}{0.4pt}%
\rule{0.4pt}{6pt}
\begin{align}
		\psi^x_i \left(x_{i,t+1}\right) &= K^{xx}_{ii} \left(\psi^x_i \left(x_{i,t}\right) - \sum_{j \neq i} K^{xx}_{ij} \psi_j \left(x_{j,t}\right) - K^{xw}_i \frac{1}{m} \sum \limits_{n=0}^{m-1} \psi^w_i \left(w_{i,t + n \tau} \right) - K^{xy}_i \frac{1}{m} \sum \limits_{n=0}^{m-1} \psi^y_i \left(y_{i,t + n \tau} \right) \right) \nonumber \\
		&+ \sum_{j \neq i} K^{xx}_{ij} \psi_j \left(x_{j,t}\right) + K^{xw}_i \frac{1}{m} \sum \limits_{n=0}^{m-1} \psi^w_i \left(w_{i,t + n \tau} \right) + K^{xy}_i \frac{1}{m} \sum \limits_{n=0}^{m-1} \psi^y_i \left(y_{i,t + n \tau} \right)  \label{eq:x dyn}\\
		\psi^y_i \left(y_{i,t + \left(n+1\right) \tau}\right) &= K^{yy}_{i} \left( \psi^y_i \left(y_{i,t + n \tau}\right) - K^{yx}_i \psi^x_i \left(x_{i,t} \right) - K^{wu}_i \psi^u_i \left(u_{i,t}\right) \right) \nonumber \\
		& + K^{yx}_i \psi^x_i \left(x_{i,t} \right) + K^{wu}_i \psi^u_i \left(u_{i,t}\right) \label{eq:y dyn} \\
		\psi^w_i \left(w_{i,t + \left(n+1\right) \tau}\right) &= K^{ww}_{i} \left( \psi^w_i \left(w_{i,t + n \tau}\right) - K^{wx}_i \psi^x_i \left(x_{i,t}\right) - K^{wy}_i \psi^y_i \left(y_{i,t + n \tau}\right) - K^{wu}_i \psi^u_i \left(u_{i,t}\right) \right) \nonumber \\
		& + K^{wx}_i \psi^x_i \left(x_{i,t}\right) + K^{wy}_i \psi^y_i \left(y_{i,t + n \tau}\right) + K^{wu}_i \psi^u_i \left(u_{i,t}\right) \label{eq:w dyn}
\end{align}
\par\noindent\rule{\dimexpr(0.5\textwidth-0.5\columnsep-0.4pt)}{0.4pt}%
\rule[-6pt]{0.4pt}{6pt}
\begin{multicols}{2}

The running cost approximation is

\begin{align}
		\hat{L}_i &= \sum_t \frac{1}{m} \sum_{n = 0}^{m-1} C_{i,t + n \tau} \\
		C_{i,t + n \tau} &= \sum_{j,k} \left[\psi^j_{i,t + n\tau}\right]^T Q_{i,jk} \psi^k_{i,t+ n \tau} \nonumber \\
		&+ \sum_j \left[c^{\left(i\right)}_j\right]^T \psi^j_{i,t + n \tau} + c_{i,0} \\
		\lim_{\epsilon \rightarrow 0} \frac{1}{m} \sum_{n=0}^{m-1} C_{i,t + n \tau} &= \sum_{j,k} \left[\psi^j_{i,t}\right]^T Q_{i,jk} \psi^j_{i,t} \nonumber \\
		&+ \sum_j c^T_{i,j} \psi^j_{i,t} + c_{i,0} \\
		\psi^j_{i,t + n \tau} & = \psi^j_{i,t} \equiv \psi^j_i \left(j_{i,t}\right) \ j \in \left\{x,u\right\} \\
		\psi^j_{i,t + n \tau} & \equiv \psi^j_i \left(j_{i,t + n \tau} \right) \ j \in \left\{w,y\right\} \\
		\psi^j_{i,t} & \equiv \psi^j_i \left(j^*_i \left(x_{i,t},u_{i,t}\right)\right)  \ j \in \left\{w,y\right\}
\end{align}

We can solve for the equilibrium of the fast dynamics to produce the $\epsilon \rightarrow 0$ slow dynamics of the form \cite{bakker2025time}

\begin{align}
		\psi^w_i \left(w_i^* \left(x_{i,t},u_{i,t}\right) \right) &= B^{wx}_i \psi^x_i \left(x_{i,t}\right) + B^{wu}_i \psi^u_i \left(u_{i,t}\right) \label{eq:comb dyn 1}\\
		\psi^y_i \left(y_i^* \left(x_{i,t},u_{i,t}\right) \right) &= B^{yx}_i \psi^x_i \left(x_{i,t}\right) + B^{yu}_i \psi^u_i \left(u_{i,t}\right) \\
		\psi^x_i \left(x_{i,t+1}\right) &= B^{xx}_i \psi^x_i \left(x_{i,t}\right) + B^{xu}_i \psi^u_i \left(u_{i,t}\right) \nonumber \\
		&+ \sum_{j \neq i} \left(I - K^{xx}_{ii} \right) K^{xx}_{ij} \psi^x_j \left(x_{j,t}\right) \label{eq:comb dyn 3}
\end{align}
\begin{align}
		\boldsymbol \psi \left(\mathbf{x}_{t+1}\right) &= B_{xx,comb} \boldsymbol \psi \left(\mathbf{x}_t\right)+ B_{xu,comb} \boldsymbol \psi \left(\mathbf{u}_t\right) \label{eq:comb dyn 4}
		%
\end{align}

\noindent for each agent $i$.  Again, we can solve for the optimum or equilibrium of this problem.  The differences between the optimality conditions associated with the two solutions are analogous to those presented in the previous section.  As in our previous paper \cite{bakker2025time}, we use the $\epsilon \rightarrow 0$ limit to collapse the dynamics into the slow scale for optimization purposes.  This significantly reduces the size and complexity of the problem, which is even more significant when solving for the equilibrium than the optimum.

\subsection{Analytical Tools}

The KO can support the use of a wide range of analytical tools for linear systems.  For a dynamical system $x_{t+1} = A x_t$, the maximum initial growth of a transient is $\log \left\| \mathbf{A}\right\|_2$ and the maximum transient growth lower bound is \cite{trefethen2005spectra}

\begin{gather}
		T_{bound} = \sup \limits_{\left|z\right| > 1} \left(\left|z\right| - 1\right) \left\| \left(zI - A\right)^{-1}\right\|_2
\end{gather}

We can compare the results obtained with, e.g., $K^{xx}_{ii}$ vs. $K_{xx,comb}$ to understand how feedbacks between agents affect overall system stability and transient growth and compare the $\epsilon \neq 0$ dynamics to the $\epsilon \rightarrow 0$ dynamics to see how cross-scale feedbacks within a given agent's dynamics affect the coupling between agents at the slow scale.  We can also use control theoretic tools even when there are no control variables present.  For example, we can quantify the ``controllability'' of each agent by to the other agents using the ``controllability'' gramian $X^c_{ij}$ calculated from:
				
\begin{gather}
		K_{ii} X^c_{ij} K^T_{ii} - X^c_{ij} = -  \left(I - K_{ii} \right) K_{ij} K^T_{ij} \left(I - K_{ii}\right)^T
\end{gather}

This gramian measures how strongly variations in the state of agent $j$ affect the behavior of agent $i$ (in terms of input and output energy).  In the case with time scale separation and hierarchical control, we could compare the gramian calculated using $K^{xx}_{ii}$ vs. $B^{xx}_i$ (the ``control'' matrix $\left(I - K^{xx}_{ii}\right) K^{xx}_{ij}$ would be the same in both cases) to see how the feedbacks between time scales change (or not, as the case may be) the impact one agent can directly have on another.  The actual controllability gramians of each agent, $X^c_i$, and of the system as a whole, can be calculated using an analogous approach to quantify controllability with respect to each agent's control inputs (and combinations thereof).  We can also quantify the stability of the whole system with respect to perturbations on each agent.  To do this, consider just the uncontrolled dynamics in, e.g., (\ref{eq:comb dyn 4}) and take the SVD of $B_{xx,comb}$.  To test the impact of agent $j$, let us assume that $\psi^x_j \neq 0$ and $\psi^x_i = 0, i \neq j$.  Then the optimization is

\begin{gather}
		P_{max,i} = \max_{\psi_j} \max_k \left\| \left[B_{xx,comb}\right]^k \boldsymbol\psi^x\right\|^2 \\
		\left\| \psi^x_j \right\| = 1, \left\| \psi^x_i \right\| = 0 \ i \neq j
\end{gather}

We can solve this by using the SVD of $\left[B_{xx,comb}\right]^k$ and evaluating the fixed $k$ solution over a sufficiently large range of $k$ values; if the system is stable, the problem is bounded.  Note, however, that this optimization is only an upper bound, as there may not be an $x_j$ value corresponding to each value of $\psi^x_{j}$.  Again, we can consider these perturbations for the slow time scale dynamics on their own as well as the $\epsilon \rightarrow 0$ dynamics to understand the implications of feedbacks between time scales.

\section{Computational Demonstration}

\subsection{Problem Formulations}

Let us define a two-agent system where the dynamics are a mixture of nonlinear oscillators (Duffing, van der Pol, and nonlinearly damped pendulum) with nonlinear couplings across scales, where $\epsilon = 100$, and agents and where the actuators are governed by Proportional-Integral (PI) control loops; for the derivations of the PI loop dynamics, see \cite{bakker2025time}.  The dynamics for the full system are

\begin{align}
		\dot{x}_{1,1} &= x_{1,2} \\
		\dot{x}_{1,2} &= -\left(1 - x_{1,1}^2\right) x_{1,2} - x_{1,1} + 0.25 \sinh x_{2,1} \nonumber \\
		&+ 0.5y_{1,1} \\
		\dot{y}_{1,1} &= \frac{1}{\epsilon} y_{1,2}  \\
		\dot{y}_{1,2} &= \frac{1}{\epsilon} \left( -2 y_{1,2} - y_{1,1} - y_{1,1}^3 - 2 w_1 + 0.5 x_{1,2}^2\right) \\
		\dot{w}_1 &= \frac{1}{\epsilon} \left[y_{1,2} + \left(y_{1,1} - u_1 \right)\right] \\
		\dot{x}_{2,1} &= x_{2,2} \\
		\dot{x}_{2,2} &= - \sinh x_{2,2} - \sin x_{2,1} + 0.5 \tanh x_{1,2} + 0.5 y_{2,1} \\
		\dot{y}_{2,1} &= \frac{1}{\epsilon} y_{2,2} \\
		\dot{y}_{2,2} &= \frac{1}{\epsilon} \left( -2 y_{2,2} - y_{2,1} - y_{2,1}^3 - 2 w_2 + 0.5 x_{2,2} ^2\right) \\
		\dot{w}_2 &= \frac{1}{\epsilon} \left[y_{2,2} + \left(y_{2,1} - u_2 \right)\right]
\end{align}

\noindent where $i_{j,k}$ indicates the the $k$-th variable of type $i$ from agent $j$.  The coupling functions between scales essentially mean that the kinetic energy of the fast scale oscillator affects the dynamics of the slow scale oscillator ($0.5 x_{i,2}^2$, the more energy, the less stable the system is), the position of the fast scale oscillator affects the dynamics of the slow scale oscillator ($0.5 y_{i,1}$, acting like a force), the actuator input acts like a force on the fast scale dynamics ($2 w_i$), the position of Agent 2's slow scale oscillator nonlinearly affects the dynamics of Agent 1's slow scale oscillator ($0.25 \sinh x_{2,1}$), and the velocity of the Agent 1's slow scale oscillator nonlinearly affects the dynamics of the Agent 2's slow time scale oscillator ($0.5 \tanh x_{1,2}$).  When we consider the case without hierarchical control or time scale separation, we substitute $u_i$ for $y_{i,1}$ in the slow dynamics -- the fixed point to which the PI control drives the fast dynamics -- and then remove all of the fast scale equations and variables (states and actuators) from the formulation.  The running cost for each agent $i$ is

\begin{gather}
		L_i = x_{i,1}^2 + x_{i,2}^2 + y_{i,1}^2 + y_{i,2}^2+ w_i^2
\end{gather}

In the case without hierarchical control or time scale separation, this becomes $x_{i,1}^2 + x_{i,2}^2 + u_i^2$.  Note that this objective separates nicely along disciplinary lines.  We also constrain the slow time scale states and control inputs to lie between -1 and 1.  In the full problem, $w$ and $x$ do not appear directly in each other's dynamics, so we can use this domain information to set $K_{xw}$ and $K_{wx}$ to be identically zero (see \cite{bakker2022deception} for more on this kind of domain knowledge incorporation).

\subsection{Koopman Implementation}

To learn our Koopman representations -- one for each case -- we used state-inclusive observables for $\psi_x$, $\psi_y$, and $\psi_w$ with the nonlinear observables calculated using small neural networks, and we specified $\hat{L}_i$ as a Riemann approximation to $L_i$ \cite{bakker2025time}.  We used 12 nonlinear observables for $\psi_x$ and $\psi_y$, and four nonlinear observables for $\psi_w$; $u$ was kept unlifted, as is common in the KO literature (e.g., \cite{korda2018linear,king2021solving}).  To train the neural networks (implemented using Keras in Tensorflow), we chose 1e4 random initial conditions in the $\left[-1,1\right]^n$ hypercube and simulated the system $\Delta t = 0.1$ seconds forward in time to produce our training data; for the slow time scale data, we only recorded data at $t = 0$ and $t = \Delta t$, but for the fast time scale data, we recorded data at $\Delta t/\epsilon$ time steps during that interval.  The training function consisted of absolute prediction error and stability enforcement terms in all three cases; the stability enforcement involved calculating the eigenvalues of the relevant matrices and penalizing eigenvalues with magnitudes greater than one, and it was applied both to all of the relevant matrices in the $\epsilon > 0$ case (e.g., $K_{xx}$, $K_{yy}$) and to the slow time scale matrices produced in the $\epsilon \rightarrow 0$ limit (e.g., $K_{comb}$, $B_{xx}$).  In the case with hierarchical control and time scale separation, it was difficult to get good predictions, so we incorporated the $\epsilon \rightarrow 0$ limit into the training prediction: the model still learned the fast and slow dynamics separately, but it calculated the $B$ matrices and used them in parallel with the $\epsilon \neq 0$ model predictions.  This seemed to help with learning the fast dynamics, too, and did not involve any extra knowledge about the system (just the KO structure).  The optimal control and game theory problems were formulated in Pyomo \cite{hart2017} and solved using IPOPT \cite{wachter2006implementation} and PATH \cite{dirkse1995path}, respectively, for 100 different randomly generated initial conditions.  For comparison, we also produced baseline results assuming constant $u=0$ control inputs.

\section{Results}

\subsection{Multi-Agent System}

We produced a highly accurate KO model of the multi-agent system: the mean RMS error over 100 trajectories of 100 time steps each was 0.038 for Agent 1 and 0.034 for Agent 2.  There are two key highlights from the analyses in Table \ref{tab:stability multi-discp}.  Firstly, the interactions between agents seem to amplify transients: transient growth metrics are greater in the combined system than they are for either agent separately.  The implication is that there is some positive feedback between the agents.  However, the combined system is more controllable (as measured by $\left\|X^c_i\right\|$) than either agent is separately.  Secondly, each agent has roughly the same impact on the overall system ($P_{max,i}$) and the other agent ($\left\|X^c_{ij}\right\|$).  Note that these metrics represent simplifications of complex phenomena, so small differences in a single metric should not be over-analyzed.  The amplification of transients is supported by multiple metrics, for example, so that conclusion is more strongly supported than, say, the conclusion that Agent 1 has more influence on Agent 2 than the other way around.

\begin{center}
\begin{table}[H]
    \centering
    \caption{Analysis Results, Multiple Agents}
    \label{tab:stability multi-discp}
    \vspace{6pt}
    \begin{tabular}{l|ccc}
        Metric &Agent 1 &Agent 2 &Combined System\\ 
        \hline
        $\left\|X^c_i\right\|$ &0.018 &0.014 &0.046 \\
				$\left\|X^c_{ij}\right\|$ &1.97 &2.32 &--\\
				$P_{max,i}$ &54.67 &52.45 &--\\
				$\log \left\| \mathbf{A}\right\|_2$ &0.429 &0.194 &0.436 \\
				$T_{bound}$ &3.604 &3.581 &4.13 \\
    \end{tabular}
\end{table}
\end{center}

The equilibrium and optimum solutions produced small, similar agent control policies: for both agents, average RMS values were ${\approx} 0.05$ and the average RMS differences between optimum and equilibrium solutions were ${<}0.005$.  The sum of the agent objectives was correspondingly similar for both solutions (${<}0.4\%$ difference in all instances).  There were slight differences, however, in how the agents fared in the equilibrium and optimum solutions.  Agent 1 did slightly better, and Agent 2 slightly worse, in the equilibrium compared to the optimum, but the differences were ${<}0.2\%$.  Both agents also did slightly better in the equilibrium and optimum solutions than in the constant $u=0$ baseline: a 0.4\% improvement for Agent 1 and a 1.4\% improvement for Agent 2.  These similarities mask two key differences between the agents, though: Agent 1's objectives were ${\sim} 50\%$ larger than Agent 2's, on average, while Agent 2's control policies were ${\sim}40\%$ larger (in average RMS values) than Agent 1's.  Agent 2 ``spent'' more on control but achieved better results.  The dynamic game was thus neither fully cooperative nor symmetric.

\subsection{Multi-Agent System with Hierarchical Control and Time Scale Separation}

Adding hierarchical control and time scale separation significantly increased system complexity: the KO model was more difficult to learn, and the prediction accuracy was noticeably worse than the previous case.  The mean RMS error over 100 trajectories of 100 time steps each for the $\epsilon \rightarrow 0$ model -- the model used for the control policy calculations -- was 0.152 for Agent 1 and 0.175 for Agent 2.  As such, the stability analysis results (Table \ref{tab:stability multi-agent}) should be taken with a larger grain of salt than in the previous case.  In Table \ref{tab:stability multi-agent}, $K^{xx}_{ii}$ indicates the $\epsilon \neq 0$ slow dynamics for Agent $i$, $K_{comb}$ indicates the combined $\epsilon \neq 0$ slow dynamics for the system as a whole, $B^{xx}_i$ indicates the $\epsilon \rightarrow 0$ dynamics for Agent $i$, and $B_{xx,comb}$ indicates the combined $\epsilon \rightarrow 0$ dynamics for the system as a whole; not all metrics apply to each set of dynamics.  

Firstly, the feedback between time scales and agents seems to magnify transients slightly: the transient growth metrics are slightly larger for $B^{xx}_i$ than $K^{xx}_{ii}$, for $B_{xx,comb}$ than $K_{comb}$, for $K_{comb}$ than $K^{xx}_{ii}$, and for $B_{xx,comb}$ than $B^{xx}_i$.  As in the previous case, the combined system is also more controllable (as measured by $\left\|X^c_i\right\|$) than either agent is separately.  Controllability is not defined at the slow scale in the $\epsilon \neq 0$ case because $u$ does not appear in the slow dynamics.  Secondly, each agent has roughly the same ability to perturb the overall system, and the effect is not strongly affected by the cross-scale feedbacks (i.e., $P_{max,i}$ is very similar for $K^{xx}_{ii}$ and $B^{xx}_i$).  Finally, it seems that Agent 1 has the ability to impact Agent 2 slightly more than the other way around regardless of cross-scale feedbacks (as indicated by $\left\|X^c_{ij}\right\|$).

\begin{center}
\begin{table}[H]
    \centering
    \caption{Analysis Results, Multiple Agents with Hierarchical Control and Time Scale Separation}
    \label{tab:stability multi-agent}
    \vspace{6pt}
    \begin{tabular}{l|cccccc}
        Metric &$K^{xx}_{11}$ &$K^{xx}_{22}$ &$K_{comb}$ &$B^{xx}_1$ & $B^{xx}_2$ &$B_{xx,comb}$ \\
        \hline
				$\left\|X^c_i\right\|$ &--&--&--&0.0117 &0.0093 &0.0295 \\
				$\left\|X^c_{ij}\right\|$ &4.529 &5.817 &--&4.043 &6.403 &-- \\
				$P_{max,i}$ &97.39 &103.85 &--&103.14 &102.12 &-- \\
				\hline
				$\log \left\| \mathbf{A}\right\|_2$ &0.681 &0.402 &0.698 &0.692 &0.415 &0.709 \\
				$T_{bound}$ &5.311 &4.421 &5.880 &5.519 &4.334 &6.108\\	
    \end{tabular}
\end{table}
\end{center}

The optimum and equilibrium calculations suffered somewhat from the reduction in model accuracy, but the trends were similar to the previous case.  The average RMS difference between the optimum and equilibrium control policies was roughly 0.02 for both agents (i.e., more different than the previous problem): both performance and behavior were similar (as before).  Again, the cost and effort split between agents was of the greatest interest.  The combination of model error and increased model complexity (due to the hierarchical control and time scale separation) accentuated two previously observed trends and altered a third.  Firstly, on average, the equilibrium and optimum objective values for Agent 2 differed by less than a percentage, but the equilibrium solution produced a 5\% average improvement for Agent 1.  Secondly, Agent 1's average objective was almost double Agent 2's for both the equilibrium and optimum solutions.  Finally, Agent 2's equilibrium and optimum control policies' average RMS values were 10-15\% smaller than Agent 1's but provided ${\sim}10\%$ improvements compared to the $u=0$ case.  Agent 1 had a 2\% performance reduction in the optimum case compared to the $u=0$ case, and the equilibrium solution produced a 3.5\% improvement compared to the $u=0$ case.  In other words, the equilibrium and optimum solution provided similar overall performance but divided up the costs and benefits somewhat differently.  Overall, in distinction to the previous problem, Agent 1 had to ``work harder'' for worse results than Agent 2.

\section{Discussion and Conclusions}

This paper presented KO model formulations for multi-agent systems with and without hierarchical control and time scale separation.  Those formulations facilitated cross-agent and cross-scale analyses.  Other, more complicated analyses like the Koopman-based causality analysis of Rupe et al. \cite{rupe2024causal} could also be employed using the Koopman formulation presented here.  Given the difficulties in learning the full $\epsilon \neq 0$ dynamics, if the information about the fast scale is not needed, it might be more expedient to learn the $\epsilon \rightarrow 0$ model directly. 

The paper also demonstrated the successful calculation of both optimal and equilibrium solutions in a general-sum game theoretic context.  These results would have been more interesting had the gap between the optimal and equilibrium solutions been greater (e.g., a more adversarial situation) and had the resulting control policies been more variable, but the methodology was sound.  The most interesting observations had to do with the relationship between agent effort and agent benefit in the equilibrium vs. optimum solutions -- the game was not symmetric, but it was almost cooperative.

There is also an important computational consideration to the methodology: the optimum solved in tenths of a second while the equilibrium could take 4 minutes to over three hours (averaging around 1-1.5 hours).  This is partially due to the fact that the MCP is always a lot larger than the optimization (due to the inclusion of dual variables), but the MCP is also a harder computational problem to solve.  There might be ways of mitigating that increased cost, though.  For example, if the equilibrium and optimum solutions are very similar, the optimum could be used as a warm start for the MCP.  Ultimately, the equilibrium calculation is only worthwhile if we expect there to be large differences between the optimum and the equilibrium solutions, but it may not be possible to know this \textit{a priori}.  

As far as game theoretic solutions go, though, the equilibrium solutions calculated and discussed here are only the tip of the iceberg, and thus there is a wealth of potential future work.  Using KO-based optimality conditions, it is possible to consider other forms of strategic interactions.  A leader-follower game, for example, can be formulated as a Mathematical Program with Equilibrium Constraints (MPEC) and solved using standard nonlinear programming techniques \cite{ruiz2014tutorial}.  MPEC-based solution methods might even be useful as a more efficient solution process in the Nash equilibrium case presented here, though we did not explore that in this paper.  Another natural extension of this work is to larger systems consisting of greater numbers of more complex agents and to the incorporation of PI control alternatives (as discussed in \cite{bakker2025time}) for lower level controls.



\bibliographystyle{ieeetr}
\bibliography{bib}

\end{multicols}

\end{document}